\documentclass[1p]{elsarticle}

\usepackage{lineno,hyperref}
\usepackage{pgf,tikz}
\usepackage{mathrsfs}
\usetikzlibrary{arrows}

\modulolinenumbers[5]

\journal{Journal of \LaTeX\ Templates}

\newtheorem{thm}{Theorem}
\newtheorem{lem}{Lemma}

\newdefinition{algoritmo}{Algorithm}
\newdefinition{df}{Definition}
\newdefinition{rmk}{Remark}
\newproof{pf}{Proof}
\newproof{pot}{Proof of Theorem \ref{thm2}}

\usepackage{verbatim}

\usepackage{amsmath}
\usepackage{pgf,tikz}
\usetikzlibrary{arrows}
\begin{document}

\begin{frontmatter}

\title{Structure of cycles in Minimal Strong Digraphs\tnoteref{mytitlenote}}
\tnotetext[mytitlenote]{This work was presented in the 10th Andalusian Meeting
on Discrete Mathematics (2017).}
%\tnotetext[mytitlenote]{Fully documented templates are available in the elsarticle package on
%\href{http://www.ctan.org/tex-archive/macros/latex/contrib/elsarticle}{CTAN}.}

%% Group authors per affiliation: %\author{Arcos-Argudo, Miguel\fnref{myfootnote}} %\address{Radarweg 29, Amsterdam}

\author[ups,upm]{Miguel Arcos-Argudo\corref{correspondingauthor}}
\cortext[correspondingauthor]{Corresponding author.}
\ead{marcos@ups.edu.ec}

\author[upm]{Jes\'us Garc\'ia-L\'opez}
\ead{jglopez@etsisi.upm.es}

\author[upm]{Luis M. Pozo-Coronado}
\ead{luispozo@etsisi.upm.es}

\address[ups]{Grupo de Investigaci\'on en Inteligencia Artificial y
Tecnolog\'{i}as de Asistencia, Universidad Polit\'ecnica Salesiana del
Ecuador}
\address[upm]{Dep. de Mat. Aplicada a las Tecnolog\'{i}as de la Informaci\'on y
las Comunicaciones, ETS de Ingenier\'{i}a de Sistemas Inform\'aticos,
Universidad Polit\'ecnica de Madrid}

\begin{abstract}
This work shows a study about the structure of the cycles contained in a
Minimal Strong Digraph (MSD). The structure of a given cycle is determined by
the strongly connected components (or strong components, SCs) that appear after
suppressing the arcs of the cycle. By this process and by the contraction of
all SCs into single vertices we obtain a Hasse diagram from the MSD.
Among other properties, we show that any SC conformed by
more than one vertex (non trivial SC) has at least one linear vertex (a vertex
with indegree and outdegree equal to 1) in the MSD (Theorem 1); that in the Hasse diagram
at least one linear vertex exists for each non trivial maximal (resp. minimal) vertex
(Theorem 2); that
if an SC contains a number $\lambda$ of vertices of the cycle then it contains at
least $\lambda$ linear vertices in the MSD (Theorem 3); and, finally, that given a cycle of
length $q$ contained in the MSD, the number $\alpha$ of linear vertices
contained in the MSD satisfies $\alpha \geq \lfloor (q+1)/2 \rfloor$ (Theorem 4).
\end{abstract}


\end{frontmatter}

%\linenumbers

\section{Introduction}

\noindent A minimal strong digraph (MSD) is a strong digraph in which the deletion of any arc yields a non strongly connected digraph. In~\cite{cita1} a compilation of properties of MSDs can be found, and in~\cite{cita2} an update of the catalog of properties is given, together with a comparative analysis between MSDs and non directed trees, from which a striking analogy between the two families of graphs emerges.

The configuration of cycles of an MSD determines, among many other things, the characteristic polynomial of its adjacency matrix and, therefore, it is very important for the spectral theory of this family of graphs~\cite{cita1,cita2}. It is also related to algorithmic problems that, restricted to MSDs, can present unexpected properties, given the strong analogy that exists between this family of graphs and trees~\cite{cita2}. In this way it is interesting to study the length $l$ of the longest cycle in an MSD, relate this parameter with the number of arcs and design an algorithm that calculates it quickly. Given the number $n$ of vertices and the number $m$ of arcs of the MSD, it is known~\cite{cita4} that
$$
\left\lceil\dfrac{m}{m-n+1}\right\rceil\leq l\leq 2n-m
$$
and that these bounds are tight. The authors also relate the parameter $l$ to the number of linear vertices~\cite{cita0} of the MSD.

The structure of a cycle $C_q$ of length $q$ in an MSD $D$ is determined by the strongly connected components (strong components, SCs) of the digraph $D^\prime$ that appear after suppressing the arcs of the cycle. The contraction of these SCs generates an acyclic digraph that, having no transitive arcs because it comes from an MSD, is in fact a Hasse diagram. We analyze the properties of this Hasse diagram and its vertices (which are SCs of the digraph $D'$).

The outline of the article is as follows: In section 2 we set up notations and discuss some basic properties. In section 3 we introduce some basic properties of the SCs of $D^\prime$ and present an algorithm that calculates the number of configurations of these SCs. In section 4 we prove the results included in the abstract. Finally, in section 5 we draw some conclusions.

\section{Notation and basic properties}

In this paper we use some concepts and basic results about graphs that are
described below, in order to fix the notation~\cite{BJG, cita1, cita2}.

Let $D=(V,A)$ be a digraph. If $(u,v)\in A$ is an arc of $D$, we say that $u$ is the initial vertex and $v$
the final vertex of the arc and we denote the arc by $uv$. We shall consider only directed paths and directed cycles.

In a strongly connected digraph, the \emph{indegree} and the \emph{outdegree} of
every vertex are greater than or equal to $1$. We shall say that $v$ is a \emph{linear vertex} if it satisfies $d^+(v)=d^-(v)=1$.

An arc $uv$ in a digraph $D$ is \emph{transitive} if there exists another $uv$-path disjoint to the arc $uv$. A digraph is
called a \emph{minimal} digraph if it has no transitive arcs.

The \emph{contraction} of a
subdigraph consists in the reduction of the subdigraph to a unique vertex. Note that the
 contraction of a cycle of length $q$ in an SD yields another SD. In such a process, $q-1$ vertices and $q$ arcs are eliminated.

A vertex $v$ on a digraph is called a \emph{cut point} if the deletion of $v$ and all its incident arcs yields a disconnected digraph.

Some basic properties concerning MDSs can be found in \cite{BJG,cita1,cita2,cita3}.
We summarize some of them: In an MSD with $n$ vertices and $m$ arcs, $n\leq m\leq 2(n-1)$; the contraction of a cycle in an MSD preserves the
minimality, that is, it produces another MSD; hence, if we contract a strongly connected subdigraph in a minimal digraph, the resulting digraph is also minimal; each MSD of order $n\geq 2$ has at least two linear vertices;
in an MSD with exactly two linear vertices, each of them belongs to a unique
cycle; furthermore, these cycles are $C_2$ or $C_3$.
The next result will be explicitly used in some of our proofs.
\begin{lem}\label{lemacut} \cite{cita2}
If an MSD contains a cycle $C_2$, then the vertices on the cycle are linear
vertices or cut points.
\end{lem}

A \emph{strong component} (SC) is a maximal strongly connected subdigraph.  A \emph{trivial strong component} is an SC that contains only one vertex. Note that every SC is a subdigraph induced by the set of its vertices.

Let $D$ be an MSD and $C_q$ a cycle contained in $D$. The \emph{digraph
associated} to $(D,C_q)$, $D'$, is the resultant digraph after suppressing in
$D$ all arcs of $C_q$. The \emph{anchored SCs} of $D'$ are the $SCs$ that
contain vertices of $C_q$. The contraction of all SCs of $D'$ yields an acyclic digraph
$H$ with no transitive arcs. Hence, $H$ is minimal and can be seen as the Hasse diagram of a partial ordering. We call $H$ the \emph{Hasse diagram associated} to $(D,C_q)$. A vertex on $H$ is \emph{minimal} (resp. \emph{maximal}) if its indegree (resp. outdegree) vanishes. Hence, minimal vertices on $H$ correspond to initial strong components of $D'$, and maximal vertices on $H$ correspond to terminal strong components of $D'$. %%%
We say that a maximal or minimal vertex of $H$ is
\emph{non trivial} if it has at least one incident arc.

\section{Basic properties of the anchored SCs of the digraph associated to
$(D,C_q)$} \label{section2}

\noindent Let $D$ be an MSD, $C_q$ a cycle contained in $D$ and $D'$ the digraph associated to $(D,C_q)$. Therefore, $D'$ is not strongly connected and its SCs may contain vertices of C, even more than one. Note that any SC is a minimal digraph.

In this section we show some results obtained through an algorithm. This algorithm has been designed to exhaustively compute all possible SC configurations of the associated digraph $D'$ for $2\leq q \leq 19$.
%Consequently the algorithm computes all possible configuration of $D'$.
We have computed only non isomorphic SC configurations of the associated digraph $D'$.

In order to compute all possible SC configurations we have considered all possible connections between vertices of $C_q$ and we have discarded the configurations in which at least one arc of $C_q$ has become transitive. Finally we have calculated and preserved only one canonical representative for each isomorphism class.

But, before passing to the description of the algorithms, we need some results on the SCs structure that we demonstrate below.

%%%%%%%%%%%%%%%%%%%%%%%%%%%%%%%%%%%%%%%%%%%%%%%%%%%%%%%%%%%%%%%%

\begin{lem}\label{lema6}
Let $D$ be an MSD, $C_q$ a cycle contained in $D$ and $D'$ the digraph associated to $(D, C_q)$. An anchored SC of $D'$ can not contain consecutive vertices of $C_q$. Hence,
the number $\lambda$ of vertices of $C_q$ in a given anchored SC of $D'$ satisfies $1\leq \lambda \leq \left\lfloor\dfrac{q}{2}\right\rfloor$.
\end{lem}
\begin{pf}
Suppose that an SC has two consecutive vertices $u$ and $v$. Then, the arc $uv$ would be transitive, contradicting the fact that $D$ is an MSD.

Let now $S$ be an SC that contains a number $\lambda$ of vertices of $C_q$. Then,
$S$ cannot contain consecutive vertices of $C_q$, and hence $\lambda \leq \left\lfloor\dfrac{q}{2}\right\rfloor$. The other inequality, $1\leq \lambda$, is straightforward.
\end{pf}

\begin{df}
Let $D$ be an MSD, $C_q$ a cycle contained in $D$ and $D'$ the digraph associated to $(D, C_q)$. We say that two anchored SCs of $D'$ are \emph{cut} if each one of them has two vertices in $C_q$, $u_1$ and $u_2$ in the first and $v_1$ and $v_2$ in the second, such that the order of these vertices in the cycle $C_q$ is $u_1,\dotsc,v_1,\dotsc,u_2,\dotsc,v_2,\dotsc,u_1$.
\end{df}

\begin{lem}\label{lema7}
Let $D$ be an MSD, $C_q$ a cycle contained in $D$ and $D'$ the digraph associated to $(D, C_q)$. If two anchored SCs of $D'$ are cut, they cannot have consecutive vertices in the cycle $C_q$.
\end{lem}
\begin{pf}
Let $u$ and $v$ be vertices contained in $C_q$. Suppose that $u$ and $v$ are consecutive vertices of $C_q$ and are contained in distinct SCs of $D'$, $u\in SC_1$ and $v\in SC_2$, that are cut. Then, the arc $uv$ would be transitive, contradicting the fact that $D$ is an MSD.
\end{pf}

%%%%%%%%%%%%%%%%%%%%%%%%%%%%%%%%%%%%%%%%%%%%%%%%%%%%%%%%%%%
Now we present the algorithms that we have implemented. Each SC configuration is represented by an integer number array with length $q$, that we have called $CompV$ and should verify the following properties:
\begin{enumerate}
\item[1)] $CompV[k]$ represents the SC to which the vertex $k$ belongs.
\item[2)] $0\leq CompV[k] < q$ for all $0 \leq k < q$ (the maximum number of components is $q$).
\item[3)] $CompV[k] \neq CompV[(k+1) \mod q]$ for all $0 \leq k < q$ (by Lemma \ref{lema6}).
\item[4)] $CompV[0] = 0$ (given a SC configuration, there is always a numbering of the SCs that verifies this).
\item[5)] $CompV[k] \leq 1+max\{CompV[j]\ |\ 0 \leq j < k\}$ for all $0 < k < q$ (given a SC configuration, there is always a numbering of the SCs that verifies this).
\end{enumerate}

The first algorithm implements the {\it ``Next"} function that, starting in the initial configuration, generates all the SC configurations. Every time it is executed it calculates the following configuration to the current one. It generates all SC configurations in lexicographical order of the $CompV$ array. The initial and final configurations are as follows:
\begin{enumerate}
\item[1)] Initial: $[0,1,0,...,0,1]$, if $q$ is even, or $[0,1,0,1,...,0,1,2]$, if $q$ is odd.
\item[2)] Final: $[0,1,2,...,q-2,q-1]$.
\end{enumerate}

\begin{algoritmo}
\emph{Next function}
\begin{enumerate}
\item[i)] Locate (looking for from ($q-1$) to $0$) the first index $k$ that allows to increase the value of $CompV[k]$. The index $k$ is the first that verifies
\begin{equation*}
CompV[k] \leq max\{ CompV[j]\ |\ 0 \leq j < k \}
\end{equation*}
\item[ii)] If $k = 0$ return NULL (There are no more configurations. The current and final configurations are the same).
\item[iii)] Assign a new value to $CompV[k]$, $comp_k$:
   \begin{equation*}
     comp_k = \left\{
	       \begin{array}{ll}
		 CompV[k]+1    & \mathrm{if\ } CompV[k-1]\neq CompV[k]+1 \\
		 CompV[k]+2    & \mathrm{if\ } CompV[k-1]= CompV[k]+1
	       \end{array}
	     \right.
   \end{equation*}
\item[iv)] Assign values to $CompV[j]$ for $k<j<q$:
   \begin{equation*}
     CompV[j] = \left\{
	       \begin{array}{ll}
		 0    & \mathrm{if\ } j\neq q-1 \mathrm{\ and\ } CompV[j-1]\neq 0 \\
		 1    & \mathrm{if\ } CompV[j-1]=0 \\
         1    & \mathrm{if\ } j=q-1 \mathrm{\ and\ } CompV[j-1]\neq 1 \\
         2    & \mathrm{if\ } j=q-1 \mathrm{\ and\ } CompV[j-1] = 1
	       \end{array}
	     \right.
   \end{equation*}
\item[v)] Return $CompV$.
\end{enumerate}
\end{algoritmo}

\begin{figure}[h]
\begin{center}
\begin{tikzpicture}[line cap=round,line join=round,>=triangle 45,x=0.34916201117318435cm,y=0.3496503496503497cm]
\clip(-4.3,-4.2) rectangle (10.02,7.3);
\draw [line width=2.pt,dash pattern=on 2pt off 2pt] (2.36,0.96) ellipse (1.7199197557352337cm and 1.7223252379110874cm);
\draw [line width=2.pt] (2.36,0.96) ellipse (0.8456931684581618cm and 0.8468759561063551cm);
\draw [line width=2.pt] (2.26,3.38)-- (1.1,4.5);
\draw [line width=2.pt] (2.12,5.88)-- (1.1,4.5);
\draw [line width=2.pt] (2.12,5.88)-- (3.28,4.68);
\draw [line width=2.pt] (3.28,4.68)-- (2.26,3.38);
\draw [line width=2.pt] (4.285675702967575,-0.5090721857690781)-- (5.84,0.04);
\draw [line width=2.pt] (5.84,0.04)-- (6.660241645845717,-1.4424824218574654);
\draw [line width=2.pt] (6.660241645845717,-1.4424824218574654)-- (5.16,-1.86);
\draw [line width=2.pt] (5.16,-1.86)-- (4.285675702967575,-0.5090721857690781);
\draw [line width=2.pt] (0.3191998793554407,-0.34442894309239946)-- (-1.28,0.02);
\draw [line width=2.pt] (-1.28,0.02)-- (-1.831694459708737,-1.6272181114948685);
\draw [line width=2.pt] (-1.831694459708737,-1.6272181114948685)-- (-0.22,-1.96);
\draw [line width=2.pt] (-0.22,-1.96)-- (0.3191998793554407,-0.34442894309239946);
\draw (3.,4.54) node[anchor=north west] {$C_4$};
\draw (-0.04,-0.96) node[anchor=north west] {$C_4$};
\draw (2.9,-1) node[anchor=north west] {$C_4$};
\begin{scriptsize}
\draw [fill=black] (2.12,5.88) circle (2.5pt);
\draw[color=black] (-2,4.7) node {$C_q$};
\draw [fill=black] (2.26,3.38) circle (2.5pt);
\draw[color=black] (1.2,2) node {$C_{2p}$};
\draw [fill=black] (0.1419012173479559,1.93285034326844) circle (2.5pt);
\draw [fill=black] (0.3191998793554407,-0.34442894309239946) circle (2.5pt);
\draw [fill=black] (2.521113382820589,-1.456700742308842) circle (2.5pt);
\draw [fill=black] (4.285675702967575,-0.5090721857690781) circle (2.5pt);
\draw [fill=black] (4.602045533471822,1.8763142615058745) circle (2.5pt);
\draw [fill=black] (1.1,4.5) circle (2.5pt);
\draw [fill=black] (3.28,4.68) circle (2.5pt);
\draw [fill=black] (5.84,0.04) circle (2.5pt);
\draw [fill=black] (6.660241645845717,-1.4424824218574654) circle (2.5pt);
\draw [fill=black] (5.16,-1.86) circle (2.5pt);
\draw [fill=black] (-1.28,0.02) circle (2.5pt);
\draw [fill=black] (-1.831694459708737,-1.6272181114948685) circle (2.5pt);
\draw [fill=black] (-0.22,-1.96) circle (2.5pt);
\end{scriptsize}
\end{tikzpicture}
\end{center}
\caption{Realization of a SC having $p$ vertices in $C_q$}
\label{RealizingSC}
\end{figure}
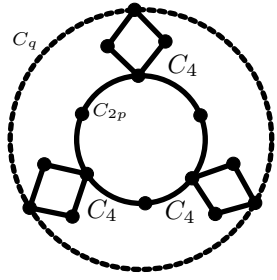

Algorithm 1 generates all the SC configurations that verify Lemma \ref{lema6}, \emph{i.e.}, configurations in which no SC contains consecutive vertices of the cycle $C_q$. According to Lemma \ref{lema7}, configurations containing SCs that are cut and have consecutive vertices of $C_q$
must be discarded. It is easy to prove that every SC configuration verifying Lemmas \ref{lema6} and \ref{lema7} is feasible as MSD. Figure~\ref{RealizingSC} shows a way to build the SCs.

To preserve only a canonical representative we have implemented the ``\emph{Canonical}'' function, that considers the $q$ arrays that are obtained when fixing each one of the positions of the circular array $CompV$ as origin. The canonical representative will be the least of them, using lexicographical order.

\begin{algoritmo}
\emph{Canonical function}
\begin{enumerate}
\item[i)] Copy $CompV$ in $CompC$ ($CompC$ is going to be the canonical array).
\item[ii)] For every $0<k<q$:\newline
$\left.\right.$\quad $\circ$ Calculate $Aux=[CompV[k],...,CompV[k+q-1]]$ (indices mod $q$).\newline
$\left.\right.$\quad $\circ$ Renumbering the SCs of $Aux$ so that $CompV$ properties hold:\newline
$\left.\right.$\quad\quad\quad $\bullet$ Initialize $ReN[j]=0$ for $0\leq j<q$ (to indicate whether vertex $j$\newline
$\left.\right.$\quad\quad\quad \,\,\, is renumbered or not).\newline
$\left.\right.$\quad\quad\quad $\bullet$ $NewComp=0$.\newline
$\left.\right.$\quad\quad\quad $\bullet$ For $j=0$ to $j=q-1$:\newline
$\left.\right.$\quad\quad\quad\quad\quad $\star$ If $ReN[j]=0$ (the vertex $j$ is not renumbered), then:\newline
$\left.\right.$\quad\quad\quad\quad\quad\quad\quad $\ast$ $Comp=Aux[j]$.\newline
$\left.\right.$\quad\quad\quad\quad\quad\quad\quad $\ast$ For $j\leq k<q$ do: if $Aux[k]=Comp$ and $ReN[k]=0$:\newline
$\left.\right.$\quad\quad\quad\quad\quad\quad\quad\quad\quad $\diamond$ $Aux[k]=NewComp$ and $ReN[k]=1$.\newline
$\left.\right.$\quad\quad\quad\quad\quad\quad\quad $\ast$ $NewComp=NewComp+1$.\newline
$\left.\right.$\quad $\circ$ If $Aux<CompC$ in lexicographical order, copy $Aux$ in $CompC$.\newline
\item[iii)] Return $CompC$.
\end{enumerate}
\end{algoritmo}

Table 1 shows the obtained results. Note how the number of SC configurations of the associated digraph $D'$ increases exponentially as $q$ grows.
\begin{center}
\begin{tabular}{rrrr}
\noalign{\hrule}
$q$ & \qquad Num. of SCs & \qquad\qquad $q$ & \qquad Num. of SCs \\
\noalign{\hrule}
2 & 1 & 11 & 162 \\
3 & 1 & 12 & 427 \\
4 & 2 & 13 & 1016 \\
5 & 2 & 14 & 2836 \\
6 & 5 & 15 & 7432 \\
7 & 6 & 16 & 20579 \\
8 & 16 & 17 & 52622 \\
9 & 28 & 18 & 159172 \\
10 & 43 & 19 & 449390 \\
\noalign{\hrule}
&
\end{tabular}

Table 1: Number of configurations of SCs of the digraph associated to $(D,C_q)$.
\end{center}

\noindent

%%%%%%%%%%%%%%%%%%%%%%%%%%%%%%%%%%%%%%%%%%%%%%%%%%%%%%%

From the experimental data, we conjecture that the number of anchored SCs of $D'$ has $\lfloor (q+3)/2\rfloor$ as a lower bound.

\section{Properties of the SCs of the digraph associated to $(D, C_q)$.}

Throughout this section, we keep using the notations introduced in the first sentence
of section~\ref{section2}, adding the Hasse diagram $H$. If $S$ is a SC of a digraph $D$,
we shall denote by $s$ the corresponding vertex of $H$. Now, we state and prove our main results.

\begin{rmk}\label{rmk1}
Let $D$ be an MSD, $C_q$ a cycle contained in $D$ and $D'$ and $H$ the digraph
and the Hasse diagram respectively associated to $(D,C_q)$. Let $s_1,\dots,s_j
$ be a path in $H$ and $u$ and $v$ vertices of $D$ such that $u \in S_1$ and $v
\in S_j$. Then it is obvious that there exists a $uv$-path in $D$ passing by
$S_1,\dots,S_{j-1} $ and $S_j$.
\end{rmk}

\begin{lem}\label{lema10}
Let $D$ be an MSD, $C_q$ a cycle contained in $D$ and $D'$ and $H$ the digraph
and the Hasse diagram respectively associated to ($D$,$C_q$). Then a cycle with
no linear vertices in $D$ can not be an SC of $D'$.
\end{lem}
\begin{pf}

First, we state that an SC $S$ corresponding
 to a maximal or minimal vertex of $H$ contains at least
one vertex of $C_q$.

In fact, if $S$ is a maximal (resp. minimal) vertex in $H$,
 there are no output (resp. input) arcs from any vertex of $S$
to any vertex of any other SC (resp. to any vertex of $S$
from any vertex of any other SC). There is at least one SC different from $S$,
because $D'$ is not strongly connected. The strong connection of $D$ implies
that there exists an arc leaving (resp. entering) $S$. This arc must belong to $C_q$
and consequently $S$ contains at least one vertex of $C_q$.

Now, let $S$ be an SC and let $u,v \in S$ be two distinct vertices such that there exist arcs $ua$ and $bv$ in $D$
with $a$,$b$ $\notin S$. Then there exists in $D$ a $uv$-path that does not
contain arcs of $S$.

Indeed, let us start in $s_1=s$ and walk through $H$, leaving $s_1$ through the arc
$s_1s_2$ (where $S_2$ is the SC containing $a$), until we find a maximal vertex on
$H$: $s_1, s_2,\dots,s_j$. According to the paragraph above, in $S_j$ there must exist at
least one vertex of $C_q$, say $w$. Then by Remark~\ref{rmk1} there exists a
$uw$-path in $D$. Now we start again in $s'_1=s$ and we walk through $H$ in opposite
direction, leaving $s$ through the arc $s'_2s'_1$ (where $S'_2$ is the SC
containing $b$), until we find a minimal vertex of $H$: $s'_k, s'_{k-1},\dots,s'_2,s'_1$.
Again,  $S'_k$ must contain at least one vertex of $C_q$, say $t$. Then we obtain a $tv$-path in $D$  (applying Remark~\ref{rmk1} to the path
$s'_k,\dots, s'_2, s'_1$). It is obvious that there exists a $wt$-path in $C_q$.
Therefore there exists a $uv$-path in $D$ (obtained by the concatenation of the
$uw$-path, the $wt$-path and the $tv$-path) that does not use arcs of $S$.

We prove now that a cycle $C$ with no linear vertices in $D$ cannot
be an SC of $D'$. Let a cycle $C$ be an SC of $D'$ and let us consider by contradiction
that it does not contain any linear vertex on $D$. In that case, each vertex of $C$
should either have an input arc from or an output arc to an SC different from $C$, or
else be contained in $C_q$. Hence $C$ must contain at least two consecutive vertices
$u$ and $v$ such that: 1) $u$ belongs to $C_q$ or there is an arc $ua$ in $D$ with $a
\notin C$; and 2) $v$ belongs to $C_q$ or there is an arc $bv$ in $D$ with $b \notin C$.

It is thus clear, as we have shown above, that in $D$ there exists a $uv$-path that does not contain arcs of $C$. As a consequence, the arc $uv$ contained in $C$ would be transitive in $D$. This fact contradicts the minimality of $D$ and completes the proof.
\end{pf}

\begin{thm}\label{thm1}
Let $D$ be an MSD, $C_q$ a cycle contained in $D$ and $D'$ and $H$ the digraph
and the Hasse diagram respectively associated to ($D$,$C_q$). Then, any non
trivial SC of $D'$ with at most one vertex of $C_q$ has at least one linear
vertex in $D$.
\end{thm}
\begin{pf}
Let $S$ be a non trivial SC of $D'$ with at most one vertex of $C_q$.
Let $k$ be the number of vertices in $S$. We prove the result by induction on $k$.

Base step: If $k=2$ then $S$ is the cycle $C_2$ and Lemma~\ref{lema10} implies the
result.

Induction step: Let $S$ be an SC conformed by $k+1$ vertices ($k \geq 2$)
which has at most one vertex of $C_q$, and suppose by induction hypothesis that
any SC with at most $k$ vertices which has at most one vertex of $C_q$ contains
at least one linear vertex in $D$.

If $S$ is the cycle $C_{k+1}$ then Lemma~\ref{lema10} implies the result.

Suppose that $S\neq C_{k+1}$. Then $S$ contains a cycle $C$ of length $p$,
 $2\leq p \leq k$. Assume that $C$ does not contain any
linear vertex in $D$ (otherwise there is nothing to prove). $S$ must
contain at least one vertex that is not contained in $C$.
The contraction of $C$ in a unique vertex, say
$z$, produces another MSD $\bar{D}$ that maintains the cycle $C_q$
(because $S$ contains at most one vertex of this cycle), and another
SC $\bar{S}$ belonging to the digraph $\bar{D'}$ associated to $(\bar{D},C_q)$.
The number of vertices of $\bar{S}$ is equal to $k-p+2$, where
$2 \leq k-p+2 \leq k$. Hence, by induction hypothesis, $\bar{S}$
must contain at least one linear vertex in $\bar{D'}$, say $u$.
Since $z$ represents $C$, which has no linear vertices, it must be
$u\neq z$ because, in $\bar{D}$,  $degree(z)\geq 3$
(if $k=2$, each vertex of $C_2$ has degree $\geq 4$, since, by lemma~\ref{lemacut},
both of them are cut points). Therefore, $u$ is a linear vertex in $D$.

The proof is thus complete.
\end{pf}

\begin{rmk}\label{rmk2}
Let $D$ be an MSD, $C_q$ a cycle contained in $D$ and $D'$ and $H$ the digraph
and the Hasse diagram respectively associated to $(D,C_q)$. Let $uv$ be an arc
of $H$. Then the corresponding SCs of these vertices $u,v$ of $H$ cannot simultaneously contain vertices of $C_q$.
\end{rmk}

%%%%%%%%%%%%%%%%%%%%%%%%%%%%%%%%%%%%%%%%%%%%%%%%%%%%%%%%%%%%%%%%

\begin{lem}\label{lema11}
Let $D$ be an MSD, $C_q$ a cycle contained in $D$ and $D'$ and $H$ the digraph
and the Hasse diagram, respectively associated to ($D$,$C_q$). Then $H$ has at
least one linear vertex  for each path joining a minimal vertex to a maximal
vertex.
\end{lem}
\begin{pf}
Let $u_0,\dotsc ,u_j$ be a $u_0u_j$-path in $H$ joining a minimal vertex $u_0$
with a maximal vertex $u_j$.

We begin by remarking that $j \geq 2$. In fact, suppose there
is an arc in $H$ joining the minimal vertex $u_0$ to the maximal
vertex $u_1$. In the proof of Lemma~\ref{lema10}
we have shown that both $U_0$ and $U_1$ must contain vertices of $C_q$,
which is impossible by Remark~\ref{rmk2}.

Let us now suppose, by contradiction, that there are no linear vertices in the
$u_0u_j$-path. Hence, $d^{-}(u_1)=1$ or else the arc $u_0u_1$ would be
transitive in $D$. In fact, if $v_1$ is the minimal vertex reached by walking in
reverse direction from $u_1$ using an arc $u'_1u_1$ different from $u_0u_1$
(such an arc exists because of $d^{-}(u_1)>1$) then a $u_0u_1$-path (not
containing the arc $u_0u_1$) can be obtained by concatenation of a
$u_0v_1$-path %%in $D-H$ (joining a maximal to a minimal)
with the $v_1u_1$-path.

Then, $d^+(u_1)>1$, since $d^-(u_1)=1$ and we are assuming that $u_1$ is not
linear. Let $u_1''\neq u_2$ be the vertex defined by the corresponding arc
$u_1u''_1\in D$.

Now the following result will be proved for all $u_i$, $2\leq i< j$:
$d^-(u_i)=1$ and there is an arc $u_iu_i^{\prime\prime}$ with
$u_i^{\prime\prime}\neq u_{i+1}$. To show this, the following reasoning is
applied iteratively for each vertex, starting from $u_2$.
 First, we remark that $d^-(u_i)=1$. Otherwise, the arc $u_{i-1}u_i$ would be
 transitive in $D$, because an $u_{i-1}u_i$-path would exist,
 not containing the arc $u_{i-1}u_i$. In fact, let $v_i$ be the minimal vertex
 that is reached walking in reverse direction from $u_i$
 starting with an arc $u_i^\prime u_i$ different from $u_{i-1}u_i$ (such an arc
 exists, since $d^-(u_i)>1$).
 Let $w_{i-1}$ be the maximal vertex reached walking from $u_{i-1}$ starting
 with the arc $u_{i-1}u_{i-1}^{\prime\prime}$.
 Then an $u_{i-1}u_i$-path would be obtained by concatenation of the
 $u_{i-1}w_{i-1}$-path with the $w_{i-1}v_{i}$-path and with
  the $v_iu_i$-path.

 Besides, $d^+(u_i)>1$ also holds, because $d^-(u_i)=1$ and, by hypothesis,
 $u_i$ is not a linear vertex.
 Let $u_i^{\prime\prime}\neq u_{i+1}$ be the vertex defined by the arc
 $u_iu_i^{\prime\prime}$ belonging to $H$.

Finally, let us show that the arc $u_{j-1}u_j$ is transitive. In fact, let
$w_{j-1}$ be the maximal vertex reached walking from $u_{j-1}$, starting with
the arc $u_{j-1}u_{j-1}^{\prime\prime}$. The $u_{j-1}u_j$-path obtained by
concatenation of the $u_{j-1}w_{j-1}$-path and the $w_{j-1}u_j$-path proves
that $u_{j-1}u_j$ is transitive. This fact contradicts the minimality of $D$
and concludes the proof.

\end{pf}
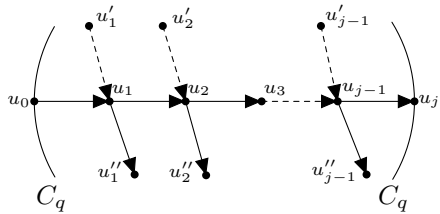
\begin{figure}[htb]
\centering
\begin{tikzpicture}[line cap=round,line join=round,>=triangle
45,x=1.0cm,y=1.0cm]
\clip(-1.91,0.17) rectangle (4.91,3.49);
\draw [->] (-0.01,2) -- (0.99,2);
\draw [->] (0.99,2) -- (1.26,1.02);
\draw [->] (0.99,2) -- (1.99,2);
\draw [->,dash pattern=on 2pt off 2pt] (1.99,2) -- (2.99,2);
\draw [->] (2.99,2) -- (4,2);
\draw [->] (-1,2) -- (-0.01,2);
\draw [shift={(0.99,2)}] plot[domain=2.62:3.67,variable=\t]({1*1.99*cos(\t
r)+0*1.99*sin(\t r)},{0*1.99*cos(\t r)+1*1.99*sin(\t r)});
\draw [shift={(1.99,2)}] plot[domain=-0.54:0.55,variable=\t]({1*2.01*cos(\t
r)+0*2.01*sin(\t r)},{0*2.01*cos(\t r)+1*2.01*sin(\t r)});
\draw [->] (-0.01,2) -- (0.32,1.03);
\draw [->,dash pattern=on 2pt off 2pt] (0.69,3) -- (0.99,2);
\draw [->,dash pattern=on 2pt off 2pt] (2.77,3) -- (2.99,2);
\draw [->] (2.99,2) -- (3.37,1.03);
\draw [->,dash pattern=on 2pt off 2pt] (-0.28,3) -- (-0.01,2);
\draw (-1.1,1) node[anchor=north west] {$C_q$};
\draw (3.4,1) node[anchor=north west] {$C_q$};
\begin{scriptsize}
\fill [color=black] (-0.01,2) circle (1.5pt);
\draw[color=black] (0.17,2.15) node {$u_1$};
\fill [color=black] (0.99,2) circle (1.5pt);
\draw[color=black] (1.17,2.15) node {$u_2$};
\fill [color=black] (1.99,2) circle (1.5pt);
\draw[color=black] (2.17,2.15) node {$u_3$};
\fill [color=black] (2.99,2) circle (1.5pt);
\draw[color=black] (3.35,2.15) node {$u_{j-1}$};
\fill [color=black] (4,2) circle (1.5pt);
\draw[color=black] (4.2,2) node {$u_j$};
\fill [color=black] (1.26,1.02) circle (1.5pt);
\draw[color=black] (0.94,1.07) node {$u''_2$};
\fill [color=black] (-1,2) circle (1.5pt);
\draw[color=black] (-1.2,2) node {$u_0$};
\fill [color=black] (-0.28,3) circle (1.5pt);
\draw[color=black] (-0.04,3.13) node {$u'_1$};
\fill [color=black] (0.32,1.03) circle (1.5pt);
\draw[color=black] (0.0,1.07) node {$u''_1$};
\fill [color=black] (0.69,3) circle (1.5pt);
\draw[color=black] (0.96,3.11) node {$u'_2$};
\fill [color=black] (2.77,3) circle (1.5pt);
\draw[color=black] (3.13,3.13) node {$u'_{j-1}$};
\fill [color=black] (3.37,1.03) circle (1.5pt);
\draw[color=black] (2.95,1.07) node {$u''_{j-1}$};
\end{scriptsize}
\end{tikzpicture}
\caption{Illustration of Lemma~\ref{lema11}.}
\end{figure}

Note that the linear vertices resulting from the lemma above cannot be the
first nor the last vertex of the original path.

Let $D$ be an MSD, $C_q$ a cycle contained in $D$, and $D'$ and $H$ the digraph
and the Hasse diagram respectively associated to $(D,C_q)$. We call
\emph{pseudominimal (resp. pseudomaximal) vertex} of $H$ to any vertex
corresponding to an anchored SC with outdegree (resp. indegree) greater than
$0$. Note that a pseudominimal vertex which is not minimal has to be also
pseudomaximal. Note also that an anchored SC must be either pseudominimal,
pseudomaximal, or trivial. Finally, we remark that the proof of Lemma~\ref{lema11} is
also valid for paths whose endpoints are pseudominimal or pseudomaximal vertices.

\begin{thm}\label{thm2}
Let $D$ be an MSD, $C_q$ a cycle contained in $D$ and $D'$ and $H$ the digraph
and the Hasse diagram respectively associated to ($D$,$C_q$). Then, the number
of linear vertices of $H$ is greater than or equal to the number of
pseudominimal (resp. pseudomaximal) vertices of $H$.
\end{thm}

\begin{pf}
Let $u_1,\dots,u_k$ be the pseudominimal vertices of $H$. For all $u_i$, $1\leq
i\leq k$, consider a $u_iv_i$-path from $u_i$ to a pseudomaximal vertex $v_i$
of $H$. Applying Lemma~\ref{lema11}, let $w_i$ be the first linear vertex that appears in
this path. By the construction used in the proof of Lemma~\ref{lema11}, the vertices
$w_1,\dots,w_k$ can be shown to be different to each other. In fact, given a
$u_i$, $1\leq i\leq k$, the vertices of the $u_iw_i$-path can not belong to any
of the paths $u_jv_j$-path ($1\leq j\leq k$, $j\neq i$) because the indegree of
each vertex is $1$ or $0$ ($u_i$). This fact implies the result in the
pseudominimal vertices case. An analog demonstration of the Lemma~\ref{lema11}, and the
previous result, but starting in the pseudomaximal vertices, concludes the proof.
\end{pf}

\begin{rmk}\label{remark3}
Lemma~\ref{lema11} also holds if we consider paths with end vertices in $C_q$ and
contained in SCs in which all cycles  include verices of $C_q$.
\end{rmk}

\begin{thm}\label{thm3}
Let $D$ be an MSD, $C_q$ a cycle contained in $D$ and $D^\prime$ the digraph
associated to $(D,C_q)$. An SC $S$ of $D'$ that contains a number $\lambda>1$
of vertices of $C_q$ has at least $\lambda$ linear vertices in $D$.
\end{thm}
\begin{pf}
Let $u_0,\dots,u_{\lambda-1}$ be the $\lambda$ vertices of $C_q$ belonging to
$S$. Let $\bar{D}$ be the MSD obtained from $D$ after iterative contraction of
the cycles that do not include  vertices of $C_q$. For each vertex $u_i$,
$0\leq i<\lambda$, consider a $u_iu_{i+1}$-path ($i$ module $\lambda$) in
$\bar{D}$. Remark~\ref{rmk3} proves that each
$u_iu_{i+1}$-path has at least a linear vertex in $\bar{D}$ and that, if $w_i$
is the first linear vertex in the $u_iu_{i+1}$-path, the indegree of all vertices of the
$u_iw_i$-path (not considering $u_i$) is equal to $1$. This implies, in turn,
that the vertices $w_0,\dots,w_{\lambda-1}$ are all different.

To get back the configuration of $D$ we can iteratively expand, in reverse
order, the vertices obtained as a result of contraction of the cycles.
Throughout this process we update the linear vertices  $w_i$, $0\leq
i<\lambda$, in the following way. Let $x$ be the vertex to expand, and $C_x$
the generated cycle as result of the expansion. If $x$ is not a linear vertex,
there is nothing to update. However, if $x$ is a linear vertex, it is necessary
to prove that a new linear vertex exists in $C_x$. Let $x_1$ and $x_2$ be the
vertices of $C_x$ that define the input arc and the output arc respectively in
$C_x$ (there are exactly $2$ because $x$ is a linear vertex). If the length of
$C_x$ is greater than or equal to $3$, or the length is $2$ and $x_1=x_2$, then
$C_x$ has a linear vertex. It remains to analyze the case where $C_x$ has
length $2$ and $x_1\neq x_2$. In this case, we have that $C_x=C_2$, neither of
the two vertices of $C_x$ is linear, because they both have a total degree of
$3$, and neither of them is a cut point because the deletion of the arcs of the
cycle does not generate two MSDs (since the vertices of $C_x$ remain with total
degree $1$). This situation contradicts Lemma~\ref{lemacut}. Therefore, this case is not
possible and the proof concludes.

\end{pf}

\begin{thm}
Let $D$ be an MSD and $C_q$ a cycle contained in $D$. Then, the number of
linear vertices of $D$ is greater than or equal to
$\left\lfloor\frac{q+1}{2}\right\rfloor$.
\end{thm}
\begin{pf}

Let $D'$ and $H$ be the digraph and the Hasse diagram, respectively, associated
to $(D,C_q)$. Let $n_1$ be the number of vertices of $C_q$ that constitute
trivial SCs, and that are not both maximal and minimal. Let $\alpha$ be the
number of linear vertices of $D$. Each one of the vertices of $C_q$ not counted
in $n_1$ belongs to a non trivial SC and, by Theorems~\ref{thm1} and \ref{thm3}, either contributes with
a linear vertex or is itself a linear vertex in $D$ (if the vertex is a trivial SC
which is both maximal and minimal). Besides, the maximum between the number of
pseudomaximals and the number of pseudominimals of $D^\prime$ is greater than
or equal to $\left\lceil\frac{n_1}{2}\right\rceil$. Therefore, by applying
Theorem~\ref{thm2}, we obtain
$$
\alpha\geq
q-n_1+\left\lceil\frac{n_1}{2}\right\rceil=q-\left\lfloor\frac{n_1}{2}\right\rfloor
$$ The minimum of the right term of the previous inequality is reached when $
n_1 = q $ and hence $$
\alpha\geq
q-\left\lfloor\frac{q}{2}\right\rfloor=\left\lceil\frac{q}{2}\right\rceil=\left\lfloor\frac{q+1}{2}\right\rfloor
$$

\end{pf}

\section{Conclusions}
The deletion of the arcs of a cycle $C_q$ in an MSD and the structure of the SCs
in the obtained digraph $D'$ gives an insight of the possible ways in which a
cycle can be imbedded into an MSD. Some interesting properties of MSDs can be
derived from this procedure. In fact, it has allowed us to prove that an MSD
containing a $q$-cycle must have at least $\lfloor (q+1)/2 \rfloor$ linear
vertices. We think that further work on this direction will lead to a deeper
understanding of the cycle structure of MSDs, and hence of the spectral
properties of this family of digraphs.

\section*{Acknowledgement}
We want to express our gratitude to the unknown referees, whose suggestions and careful reading of our manuscript have greatly improved the final result.

\section*{References}

%\bibliography{mybibfile}

\end{document}